\documentclass[preprint,sort&compress]{elsarticle}
\usepackage{amsmath}
\usepackage{amssymb}
\usepackage{amsthm}
\usepackage{xypic}
\usepackage{graphicx}
\usepackage{mathrsfs}
\usepackage{ifpdf}
\usepackage{paralist}
\allowdisplaybreaks[4]

\numberwithin{equation}{section}

\usepackage{color}
\usepackage{bm}

\ifx\red\undefined
\newcommand{\red}{\color{red}}

\renewcommand{\le}{\leqslant}
\renewcommand{\ge}{\geqslant}

\newcommand{\bR}{{\mathbb R}}

\newcommand{\mm}{\hspace{2em}}

\newcommand\beaa{\begin{eqnarray*}}
\newcommand\eeaa{\end{eqnarray*}}
\newcommand{\re}[1]{\mbox{$($\ref{#1}$)$}}

\newtheorem{thm}{Theorem}[section]

\newtheorem{cor}[thm]{Corollary}

\newcommand\bea{\begin{eqnarray}}
\newcommand\eea{\end{eqnarray}}
\newcommand\be{\begin{equation}}
\newcommand\ee{\end{equation}}
\usepackage[sc]{mathpazo}

\begin{document}

\begin{frontmatter}
\title{The existence of periodic solution and asymptotic behavior of solutions
 for a multi-layer tumor model with a periodic provision of external nutrients}

\author[sysu]{Wenhua He}
\ead{hewh27@mail2.sysu.edu.cn}
\author[sysu]{Ruixiang Xing\corref{cor}}
\ead{xingrx@mail.sysu.edu.cn}

\cortext[cor]{Corresponding author}

\address[sysu]{School of Mathematics, Sun Yat-sen University, Guangzhou 510275, China}

\begin{abstract}
In this paper, we consider a multi-layer tumor model with a periodic provision of external nutrients. The domain occupied by tumor has a different shape (flat shape) than spherical shape which has been studied widely.
The important parameters are periodic external nutrients $\Phi(t)$ and threshold concentration for proliferation $\widetilde{\sigma}$.
In this paper, we give a complete classification about $\Phi(t)$ and $\widetilde{\sigma}$ according to global stability of zero equilibrium solution or global stability of the positive periodic solution. Precisely, if  $\frac{1}{T} \int_{0}^{T} \Phi(t)d t\leqslant\widetilde{\sigma}$, then the zero equilibrium solution is globally stable while if  $\frac{1}{T} \int_{0}^{T} \Phi(t)d t>\widetilde{\sigma}$, then there exists a unique positive T-periodic solution and it is globally stable.
\end{abstract}

\end{frontmatter}

\section{Introduction}\label{sec:intro}

 In this paper, we consider a free-boundary multi-layer tumor model with a periodic provision of external nutrients. Let
$$
\Omega(t) \triangleq \{(x,y)\in \bR^{2}\times \bR; \;\;  0<y< \rho(t,x)\}, \mm {\bm x} = (x,y)= (x_1,x_2, y),
$$
be the flat-shaped region of the 3-dimensional tumor where the positive function $\rho(t,x)$ is  unknown and it is the free boundary.
Denote $\Gamma(t)$ by  the upper free boundary
$\{y= \rho(t,x)\}$ of $\Omega(t)$, which is a permeable layer and  $\Gamma_0$ by the lower boundary $\{y=0\}$, which is fixed and impermeable.
Such flat-shaped domain of tumors is used to model multi-layer tumor growth developed by medico-biologists.
It is important work to study the metabolism of multi-layer tumor tissues. For more details, we refer to the papers \cite{1997Three,2004Three,kyle1999characterization}. 

The concentration $\sigma$ of nutrient (mostly oxygen or glucose ) satisfies the reaction diffusion equation:
\begin{eqnarray}
\lambda\sigma_t - \Delta\sigma + \sigma = 0, \hspace{2em} (x_{1}, x_{2}, y)\in\Omega(t),\hspace{2em} t>0,\label{1.1}
\end{eqnarray}
with the boundary condition
\begin{eqnarray}
\displaystyle\frac{\partial \sigma}{\partial y}\Big|_{\Gamma_0} =0 , \hspace{2em} \sigma \Big|_{\Gamma(t)} = \Phi(t), \hspace{2em} t>0,\label{1.2}
\end{eqnarray}
where the positive periodic function $\Phi(t)$ is the external nutrient concentration with period T. In this paper, we assume that external nutrient concentration is a periodic function instead of constant external nutrient concentration \cite{CE1}. It is more reasonable. Here, $\lambda$ denotes the ratio of the  nutrients diffusion rate 
to the
cell proliferation rate
and $\lambda\ll1$  (see \cite{Byrne1995}). Then in this paper we consider the quasi-steady state approximation, i.e.,
$\lambda=0$.

Let $p$ be the pressure, $S$ represent the proliferation rate and $\vec{V} $ denote the velocity of tumor cell movement.
 Assuming the tumor is porous medium type and combining with the Darcy's law ($\vec{V}=-\nabla p$) and the conversation of mass ($\mbox{div} \vec V = S$), we get $ -\Delta p = S$.  Suppose $S=  \mu(\sigma-\tilde{\sigma})$, where $\mu$
is the tumor aggressiveness constant and $\tilde \sigma$
represents the threshold concentration for proliferation.
Then $p$ satisfies the following equation:
\begin{equation}\label{1.5b}
-\Delta p =\mu(\sigma-\tilde{\sigma}),\hspace{2em} (x_{1}, x_{2}, y)\in\Omega(t),\hspace{2em} t>0,
\end{equation}
with the boundary condition
\begin{eqnarray}\label{1.6}
 \displaystyle\frac{\partial p}{\partial y} \Big|_{\Gamma_0}=0, \hspace{2em} p \Big|_{\Gamma(t)}= \kappa ,\hspace{2em} t>0,
\end{eqnarray}
where $\kappa$ is the mean curvature.
Supposing the velocity field is continuous to the free boundary, then the normal velocity on
  $\Gamma(t)$ is
\begin{equation}\label{1.8}
V_n = -\nabla p\cdot n = -\frac{\partial p}{\partial n},  \hspace{2em} (x_{1}, x_{2},y)\in \Gamma(t),\hspace{2em} t>0,
\end{equation}
where $n$ is the unit outside normal vector.
The initial value of domain $\Omega(t)  $ is
$
 \Omega_0.
$

For the multi-layer tumor model with \textit{constant} external nutrient concentration, there are some interesting results. For the quasi-steady state approximation,
  Cui and Escher \cite{CE1} have established local well-posedness by means of analytic
 semigroup theory. Also under the assumption $\Phi>\widetilde{\sigma}$, they have shown that the existence and uniqueness of the positive flat stationary solution and its asymptotic behavior under non-flat perturbations.
For 2-dimensional tumor, Zhou, Escher and Cui \cite{2008Bifurcation} have studied the bifurcation of the corresponding steady problem.
In the presence of inhibitors, Zhou, Wu and Cui \cite{zwc} have derived the local existence and asymptotic behavior of flat stationary solutions under non-flat perturbations. For 2-dimensional tumor, Lu and Hu \cite{2020Bifurcation} studied the bifurcation of  tumor growth with ECM
and MDE interactions. Recently, for  tumor model with time delay, He, Xing and Hu considered the linear stability of the positive flat stationary solution under non-flat perturbations for quasi-steady state approximation in \cite{hxh1} and for general case with $\lambda\neq 0$ in \cite{hxh2}, respectively.

For the classical tumor growth models with sphere-shaped domain and a periodic external nutrients, Bai and Xu \cite{Bai2013} have shown that the zero equilibrium solution is globally stable if  $\widetilde{\sigma}>\frac{1}{T} \int_{0}^{T} \Phi(t) d t$ and if the zero equilibrium solution is globally stable, then $\widetilde{\sigma}\ge\frac{1}{T} \int_{0}^{T} \Phi(t) d t$. Also, they have proved the existence, uniqueness and stability of the positive periodic solution under the assumption $\min _{0 \le t \le T} \Phi(t)>\widetilde{\sigma} $.
For 2-dimensional problem, Huang, Zhang and Hu \cite{Huang2019} have described the linear stability of the positive T-periodic solution under non-radial perturbations. Recently, He and Xing \cite{he} have filled the gap in \cite{Bai2013}, given a complete classification about $ \Phi(t)$ and $\widetilde{\sigma}$ according to the global stability of zero equilibrium and the existence of periodic solutions, and
shown the linear stability of the positive T-periodic solution under non-radial perturbations for 3-dimensional case.

The corresponding flat problem of \re{1.1}--\re{1.8} is
\begin{align}\label{1.10a}
&\frac{\partial^{2} \sigma}{\partial y^{2}}=\sigma    && y \in(0, \rho(t)),  t>0,\\\label{1.11a}
&\displaystyle\frac{\partial \sigma}{\partial y} (0,t)=0,  \mm   \sigma(\rho(t),t)=\Phi(t), &&     t>0,\\\label{1.13a}
&-\frac{\partial^{2} p}{\partial y^{2}}=\mu(\sigma-\widetilde{\sigma})  && y \in(0, \rho(t)),  t>0,\\\label{1.14a}
&\displaystyle\frac{\partial p}{\partial y} (0,t)=0, \mm p(\rho(t),t)= 0,  &&    t>0,\\\label{1.15a}
&\frac{d \rho}{d t}=-\frac{\partial p}{\partial y}  && y=\rho(t),  t>0,\\\label{1.17a}
&\rho(0)=\rho_{0}.
\end{align}
The solution of \re{1.10a}--\re{1.11a} satisfies
\begin{align}\label{1.13}
\sigma(y, t)=\Phi(t)\frac{\cosh y}{\cosh(\rho(t))}.
\end{align}
And the solution of \re{1.13a} -- \re{1.14a} is
\begin{equation}\label{1.14}
p(y,t)=  \frac{1}{2}\mu \tilde{\sigma} y^2 +\mu\Phi(t)- \frac{1}{2}\mu\tilde{\sigma}(\rho(t))^2-\mu \Phi(t)\frac{\cosh y}{\cosh (\rho(t))}.
\end{equation}
From \re{1.15a} and \re{1.14}, system \re{1.10a}--\re{1.17a} is reduced to the following system:
\begin{align}\label{1.15}
&\frac{d \rho}{d t}(t)=\mu \rho(t)\left[\Phi(t) \frac{ \tanh( \rho(t))}{\rho(t)}-\widetilde{\sigma}\right],\\ \label{1.155}
&\rho(0)=\rho_{0}.
\end{align}
 Notice that $\rho$ is a solution of \re{1.15}--\re{1.155} if and only if  $(\sigma, p, \rho)$  is a solution of \re{1.10a}--\re{1.17a}, where $\sigma$ and $ p$ are given in \re{1.13} and \re{1.14}, respectively.

Denote
\begin{align}\nonumber
\overline{\Phi}=\frac{1}{T} \int_{0}^{T} \Phi(t) d t, \quad \Phi^{*}=\max _{0 \le t \le T} \Phi(t), \quad \Phi_{*}=\min _{0 \le t \le T} \Phi(t).
\end{align}

Our main results are the following theorems:
\begin{thm}\label{thm:1.1a'}
For any initial value $\rho_{0}>0$, \re{1.15}--\re{1.155} has a unique positive global solution $\rho$.
\end{thm}

\begin{thm}\label{thm:1.1'}
 $\widetilde{\sigma}\ge\overline{\Phi}$ if and only if the zero solution of system \re{1.15}--\re{1.155} is globally stable.
\end{thm}

\begin{thm}\label{thm:1.2}
If $\widetilde{\sigma}<\overline{\Phi}$, then the following conclusions hold: \\
(i) \re{1.15}--\re{1.155} has a unique positive T-periodic  solution $\rho_{\ast}$. \\
(ii)  There exist $\delta>0$ and $C>0$ such that for any the positive solution $\rho$,
\begin{align}
|\rho(t)-\rho_{\ast}(t)|\le C e^{-\delta t}\qquad for ~    t>0. \ \label{1.17}
\end{align}
\end{thm}


From the view of model, Theorem \ref{thm:1.1'} reflects that if the mean value of external nutrient isn't sufficient for tumor cell proliferation, then all flat-shaped tumors disappear while Theorem \ref{thm:1.2} shows that if external nutrient is sufficient, then all flat-shaped tumors grow into a T-periodic state.

\vskip 2mm

From Theorem \ref{thm:1.2}, we get  the following corollary.
\begin{cor}\label{cor:1.3333}
 If $\widetilde{\sigma}<\overline{\Phi}$, then there exists a unique positive T-periodic solution $(\sigma_{*}, p_{*}, \rho_{*})$ of \re{1.10a}--\re{1.17a} satisfies
\begin{align}\nonumber
&\sigma_{\ast}(y, t)=\Phi(t)\frac{\cosh y}{\cosh(\rho_{*}(t))}, \\\label{1.1992}
&p_{\ast}(y, t)=\frac{1}{2}\mu \tilde{\sigma} y^2 +\mu\Phi(t)- \frac{1}{2}\mu\tilde{\sigma}(\rho_{*}(t))^2-\mu \Phi(t)\frac{\cosh y}{\cosh (\rho_{*}(t))},
\end{align}
where $\rho_{*}$ is the unique positive T-periodic solution of \re{1.15}--\re{1.155}.
\end{cor}

Different from the sphere-shaped model in \cite{Bai2013} and \cite{he}, our model has flat-shaped domain which causes the various distinct computations and estimates in order to show the proofs.

The paper is organized as follows. In Section 2, we show the existence and uniqueness of the positive global solution of \re{1.15}--\re{1.155}. Also, we give the necessary and sufficient condition for global stability of zero equilibrium solution. We show the existence, uniqueness and stability of the positive periodic solution in Section 3.

\section{Proofs of Theorem \ref{thm:1.1a'} and Theorem \ref{thm:1.1'}
}
In this section, we shall give the existence and uniqueness of the positive global solution of \re{1.15}--\re{1.155} and the necessary and sufficient condition for global stability of zero equilibrium solution.
\vskip 3mm\noindent
{\it\bf Proof of Theorem \ref{thm:1.1a'}.}
 From the ODE theory, the local existence and uniqueness of the solution of \re{1.15}--\re{1.155} are obvious. Since $\frac{ \tanh \rho}{\rho}$ is strictly decreasing in $\rho$ and $\frac{ \tanh \rho}{\rho}\in(0,1)$, we have
$$
-\tilde{\sigma}\rho(t) \leq \frac{d \rho}{d t}=\rho(t)\left[\Phi(t)  \frac{ \tanh( \rho(t))}{\rho(t)}-{\tilde{\sigma}}\right] \leq ( \Phi^{*}-\tilde{\sigma})\rho(t),
$$
which implies
\begin{align}\nonumber
 \rho(0) e^{-\tilde{\sigma} t} \leq \rho(t) \leq \rho(0) e^{( \Phi^{*}-\tilde{\sigma})t}.
\end{align}
 Hence, the solution doesn't blow up or disappear at a finite time. 
 Then we get the results.
\hfill$\square$


\vskip 3mm

\noindent
{\it\bf Proof of Theorem \ref{thm:1.1'}.}
At first,  we prove the solution $\rho(t)\equiv0$ of \re{1.15} is globally stable if  $\overline{\Phi}<\tilde{\sigma}$.

 For any $t \in[0, T]$, from \re{1.15}, we have
$$
\rho(t+n T)=\rho(0) e^{\mu\displaystyle\int_{0}^{t+n T} \Phi(s) \frac{ \tanh( \rho(s))}{\rho(s)} d s-\mu\tilde{\sigma} (t+n T)}.
$$
Since $\frac{ \tanh \rho}{\rho}$ is monotone decreasing in $\rho>0$ and $\frac{ \tanh \rho}{\rho}<1$, for $\overline{\Phi}<\tilde{\sigma}$, we obtain
$$
\rho(t+n T) \le \rho(0) e^{\mu\displaystyle\int_{0}^{t+n T} \Phi(s)  d s-\mu\tilde{\sigma} n T} \le \rho(0)e^{\mu  T\overline{\Phi}}e^{\mu n T(\overline{\Phi}-\tilde{\sigma})} \rightarrow 0, \mm n \rightarrow \infty.
$$

 Next, we prove the zero equilibrium solution  of \re{1.15} is globally stable if  $\widetilde{\sigma}=\overline{\Phi}$. We shall show
\begin{align}
&\rho(t+T)\le \rho(t) \qquad for ~t>0, \label{3.33}\\
&\rho(t)\le \rho(a)e^{ \displaystyle\mu(\Phi^{*}-\widetilde{\sigma})T}  \qquad for ~  t\in[a,a+T], \qquad a\ge 0, \ \label{3.2}\\
&\liminf _{t \rightarrow+\infty}  \rho(t)=0, \ \label{3.9}\\
&\lim _{t \rightarrow+\infty} \rho(t)=0. \ \label{3.16}
\end{align}

Indeed, since $0<\frac{ \tanh \rho}{\rho}<1$, \re{1.15} implies
\begin{align}
 \frac{d \rho}{d t}\le \mu \rho(t)\left[\Phi(t)-\widetilde{\sigma}\right]. \ \label{3.222}
\end{align}
Then
$
\rho(t+T) \leqslant \rho(t)e^{ \int_{t}^{t+T}\mu\left[\Phi(t)-\widetilde{\sigma}\right] d t}=\rho(t)e^{  \displaystyle\mu\left[ \overline{\Phi}T  -\widetilde{\sigma}T\right]}=\rho(t),
$
i.e., \re{3.33} is true.

Applying \re{3.222}, we obtain
$
 \frac{d \rho}{d t}\le \mu \rho(t)\left[\Phi^{\ast}-\widetilde{\sigma}\right],
$
which implies
$$
 \rho(t)\le \rho(a)e^{ \displaystyle \mu(\Phi^{*}-\widetilde{\sigma} )(t-a)}\le \rho(a)e^{ \displaystyle \mu(\Phi^{*}-\widetilde{\sigma} )T},
$$
for  $t\in[a,a+T]$. Then \re{3.2} holds.

 We use the contradiction method to prove  \re{3.9}. Assume  that
$
\liminf _{t \rightarrow+\infty}  \rho(t)=\alpha>0.
$
For $\forall$ $\varepsilon\in(0, \alpha)>0$, there exists $M>0$ such that
\begin{align}
\rho(t)>\alpha-\varepsilon  \qquad for ~    t>M. \ \label{3.188}
\end{align}
Using \re{1.15} and the monotonicity of  $\frac{ \tanh \rho}{\rho}$ in $\rho$, we get
\begin{align}\nonumber
 \frac{d \rho}{d t}
 \le \mu \rho(t)\left[\Phi(t) \frac{ \tanh( \alpha-\varepsilon)}{\alpha-\varepsilon} -\widetilde{\sigma}\right] \qquad for ~   t>M.
\end{align}
Then
\begin{equation}\label{3.13}
\rho(t^{\ast}+nT) \le \rho(t^{\ast})e^{ \displaystyle \int_{t^{\ast}}^{t^{\ast}+nT}\mu\left[\Phi(t) \frac{ \tanh( \alpha-\varepsilon)}{\alpha-\varepsilon} -\widetilde{\sigma}\right] d t}=\rho(t^{\ast})e^{ \displaystyle\mu nT \left[\frac{ \tanh( \alpha-\varepsilon)}{\alpha-\varepsilon}\overline{\Phi}-\widetilde{\sigma}\right]},
\end{equation}
for   $t^{\ast}>M$ and  $n\geqslant1$ is an integer.
Applying \eqref{3.13} and the fact
$
\frac{ \tanh( \alpha-\varepsilon)}{\alpha-\varepsilon}\overline{\Phi}-\widetilde{\sigma}<\overline{\Phi}-\widetilde{\sigma}=0,
$
we get
\begin{align}\nonumber
 \rho(t^{\ast}+nT) \rightarrow0, \qquad n\rightarrow\infty.
\end{align}
It contracts with \re{3.188}, i.e., \re{3.9} is true.

Now we turn to \re{3.16}.
Using \re{3.9}, for $\forall$ $\varepsilon>0$, there exist $M_{0}>0$ and a sequence $t_{n}\rightarrow\infty$  such that
$
\rho(t_{n})<\varepsilon$  for     $t_{n}>M_{0}$.
\re{3.33} implies
$
\rho(t_{N}+kT)\le \rho(t_{N})< \varepsilon$
for $t_{N}>M_{0}$ and the integer $k\geqslant1$. For $t>t_{N}$, there exists $k_{0}$ such that $t\in[ t_{N}+k_{0}T,t_{N}+(k_{0}+1)T)$. Applying \re{3.2}, we have
\begin{align}\begin{array}{ll}\nonumber
\rho(t) \le \rho(t_{N}+k_{0}T)e^{ \displaystyle\mu(\Phi^{*}-\widetilde{\sigma})T} \le \varepsilon  \displaystyle e^{  \displaystyle\mu (\Phi^{*}-\widetilde{\sigma}) T}.
\end{array}
\end{align}
Then \re{3.16} holds.

As so far, we have shown that if $\widetilde{\sigma}\ge\overline{\Phi}$, $\rho(t)\equiv0$  is globally stable.
Finally, we prove that if $\rho(t)\equiv0$ of \re{1.15}--\re{1.155} is globally stable, then $\widetilde{\sigma}\ge\overline{\Phi}$.

Since  $\lim _{t \rightarrow \infty} \rho(t)=0$, for $\forall$ $\varepsilon>0$, there exists $M_{1}>0$ such that $\rho(t)<\varepsilon$ for $t \ge M_{1}$. Then the fact that $\frac{ \tanh \rho}{\rho}$ is strictly decreasing implies

\begin{align}\nonumber
\frac{d \rho}{\rho} 
\ge \mu\left[\Phi(t)  \frac{ \tanh( \varepsilon)}{\varepsilon}-\tilde{\sigma}\right], \quad t \ge M_{1}.
\end{align}
Hence
$
\frac{\rho(t+T)}{\rho(t)} \ge e^{\mu T\left(\overline{\Phi}  \frac{ \tanh( \varepsilon)}{\varepsilon}-\tilde{\sigma}\right)}.
$
 If $\overline{\Phi}>\tilde{\sigma}$, we choose $\varepsilon$  such that $ \frac{ \tanh( \varepsilon)}{\varepsilon}>\frac{\tilde{\sigma}}{ \overline{\Phi}}$.  Then
$
\frac{\rho(t+T)}{\rho(t)} >1,
$
 which contradicts to the $\lim _{t \rightarrow \infty} \rho(t)=0$. Thus $\widetilde{\sigma}\ge\overline{\Phi}$ holds.
\hfill$\square$

\section{Existence, Uniqueness and Stability of the Periodic Solution
}\label{result}

In this section, we shall prove the existence, uniqueness and stability of the periodic solution of \re{1.15}--\re{1.155}.


\vskip 3mm\noindent
{\it\bf Proof of Theorem \ref{thm:1.2}.}
Because $\widetilde{\sigma}<\overline{\Phi}\le\Phi^{\ast}$ , $0<\frac{ \tanh \rho}{\rho}<1$ and $\frac{ \tanh \rho}{\rho}$ is strictly decreasing, we have $x_{2}=(\frac{ \tanh \rho}{\rho})^{-1}(\frac{\widetilde{\sigma}}{\Phi^{\ast}})$ and $\overline{x}=\frac{(\frac{ \tanh \rho}{\rho})^{-1}(\frac{\widetilde{\sigma}}{\overline{\Phi}})}{e^{\mu(\Phi^{*}-\widetilde{\sigma})T}}$ are well defined, and
$
\overline{x}<x_{2}.
$
For each $\rho_{0}\in[\overline{x}, x_{2}]$, we define the map $F$: $[\overline{x}, x_{2}]\rightarrow \mathbb{R}$ by
 $F(\rho_{0})=\rho(T)$, where $\rho$ is the solution of \re{1.15}--\re{1.155}.

At first, we prove that $F$ maps $[\overline{x}, x_{2}]$ into $[\overline{x}, x_{2}]$. Let $\rho(0)\in[\overline{x}, x_{2}]$.
Notice that $x_{2}$ is the upper solution of \re{1.15}--\re{1.155}. Applying the comparison theorem, we have
$
\rho(t)\le x_{2} $  for $t>0.
$
Hence
\begin{align}
\rho(T)\le x_{2}. \ \label{4.3}
\end{align}
We define $\overline{\rho}$ by the solution of \re{1.15} with $\overline{\rho}(0)=\overline{x}$.
The comparison theorem implies
\begin{align}
\rho(t)\ge \overline{\rho}(t) \qquad for ~ t>0. \ \label{4.6}
\end{align}
Applying  $0<\frac{ \tanh \rho}{\rho}<1$,  we have
$
  \frac{d \overline{\rho}}{d t}
  \le\mu \overline{\rho}(t)\left[\Phi^{\ast}- \widetilde{\sigma}\right].
$
Then
\begin{equation}\nonumber
\overline{\rho}(t) \leqslant \overline{x} e^{ \displaystyle\mu(\Phi^{*}-\widetilde{\sigma}) t } \leqslant \overline{x} e^{ \displaystyle\mu(\Phi^{*}-\widetilde{\sigma}) T}=\left(\frac{ \tanh \rho}{\rho}\right)^{-1}\left(\frac{\widetilde{\sigma}}{ \overline{\Phi}}\right) \qquad for ~t\in[0,T].
\end{equation}
Applying the fact that  $\frac{ \tanh \rho}{\rho}$ is strictly decreasing, we get
$
 \displaystyle\frac{ \tanh \overline{\rho}(t)}{\overline{\rho}(t)}\geqslant \frac{\widetilde{\sigma}}{\overline{\Phi}}$ for $t\in[0, T].
$
Hence,
\begin{align*}
  \frac{d \overline{\rho}}{d t}\ge  \mu \overline{\rho}(t)\left[ \Phi(t)\frac{\widetilde{\sigma}}{\overline{\Phi}}- \widetilde{\sigma}\right] \qquad for ~ t\in(0, T),
\end{align*}
which implies
\begin{align}
\overline{\rho}(T)\ge \overline{\rho}(0)e^{\displaystyle \int_{0}^{T}\mu\left[\Phi(t)\frac{\widetilde{\sigma}}{\overline{\Phi}}-\widetilde{\sigma}\right]dt}=\overline{\rho}(0)=\overline{x}.
\label{4.10}
\end{align}
From \re{4.3}, \re{4.6} and \re{4.10}, we have
$\rho(T)\in[\overline{x}, x_{2}].$
Hence  $F$ maps $[\overline{x}, x_{2}]$ into $[\overline{x}, x_{2}]$. From the continuous dependence of the solution $\rho$ on the initial value $\rho_{0}$, we have  $F$ is continuous. Using Brouwer's fixed point theorem,  it follows  that $F$ has a fixed point $\rho_{\ast}(0)$. Then the solution $\rho_{\ast}$ of \re{1.15}--\re{1.155}  with $\rho(0)=\rho_{\ast}(0) $ is a T-periodic  positive solution. We shall prove the uniqueness of the periodic solution later.

Next we turn to  prove $(ii)$.
Let
\begin{align}\label{4.7}
\rho_{min}=\min\limits_{t>0}\left\{\rho_{\ast}(t)\right\}\qquad \text{and} \qquad \rho_{max}=\max\limits_{t>0}\left\{\rho_{\ast}(t)\right\}.
\end{align}
The uniqueness of the solution of \re{1.15}--\re{1.155} implies that $\rho_{min}>0$ and $\rho_{max}>0$.
Assume  $\rho(t)$ is the solution of \re{1.15}--\re{1.155}  with  $\rho(0)>0$. Let
\begin{align}\label{qq}
\rho(t)=\rho_{\ast}(t) e^{y(t)}.
\end{align}
In order to prove $
|\rho(t)-\rho_{\ast}(t)|\le C e^{-\delta t}$,
 we need to  prove that there exist $\delta>0$ and $C>0$ such that
\begin{align}\nonumber
|e^{y(t)}-1|\le C e^{-\delta t}\qquad for ~    t>0.
\end{align}
Taking \re{qq} into \re{1.15}--\re{1.155}, we get
\begin{align}\begin{array}{ll}
y'(t) =  \mu\Phi(t)\Big[ \displaystyle\frac{ \tanh (\rho_{\ast}(t)e^{y(t)})}{\rho_{\ast}(t)e^{y(t)}}  -\displaystyle\frac{ \tanh \rho_{\ast}(t)}{\rho_{\ast}(t)} \Big].
\end{array}
\label{4.12}
\end{align}
The uniqueness of the solution of \re{1.15}--\re{1.155} implies that $\rho(t)>\rho_{*}(t)$  if $\rho(0)>\rho _{*}(0)$ and $\rho(t)<\rho_{*}(t)$ if $\rho(0)<\rho _{*}(0)$. Then $y(t)>0$ if  $y(0)>0$ and $y(t)<0$ if  $y(t)<0$. Therefore, according to the sign of $y(t)$, we divide the arguments into two cases.

Case 1: $y(t)>0$.

Applying \re{4.12} and the strictly monotonicity  of  $\frac{ \tanh \rho}{\rho}$ in $\rho$, we have $y'(t)<0$.
Using  \re{4.12}  and the mean value theorem, we have
\begin{align}\nonumber
y'(t) e^{y(t)}=  \mu\Phi(t) (\displaystyle\frac{ \tanh \rho}{\rho})'(\zeta(t))\rho_{\ast}(t)(e^{y(t)}-1)e^{y(t)}   \le -\mu\Phi_{\ast}M_{min} \rho_{min}   (e^{y(t)}-1)\nonumber,
\end{align}
here,  $\zeta(t)\in[\rho_{\ast}(t), \rho_{\ast}(t)e^{y(t)}]\subseteq [\rho_{min}, \rho_{max}e^{y(0)}]$ 
and $M_{min}=\min\limits_{x\in [\rho_{min}, \rho_{max}e^{y(0)}]}\left\{-(\frac{ \tanh \rho}{\rho})'\right\}>0$.
Therefore,
\begin{align}
\frac{(e^{y(t)}-1)'}{(e^{y(t)}-1)}\le -\mu\Phi_{\ast}M_{min} \rho_{min} \label{4.1511}.
\end{align}
Integrating \re{4.1511} over $[0,t]$, we get
\begin{align}
e^{y(t)}-1\le (e^{y(0)}-1)e^{-\mu\Phi_{\ast}M_{min} \rho_{min} t}\qquad for ~    t>0\label{4.15111}.
\end{align}

Case 2: $y(t)<0$.

Using \re{4.12} and the strictly monotonicity  of  $\frac{ \tanh \rho}{\rho}$ in $\rho$, we obtain $y'(t)>0$. From  the mean value theorem, we get
\begin{align}\nonumber
-y'(t) e^{y(t)}=  -\mu\Phi(t)  (\displaystyle\frac{ \tanh \rho}{\rho})'(\eta(t))  \rho_{\ast}(t)(e^{y(t)}-1)e^{y(t)} \le  -\mu\Phi_{\ast} \overline{M}_{min}\rho_{min}   (1-e^{y(t)})e^{y(0)}\nonumber,
\end{align}
here,  $\eta(t)\in[\rho_{\ast}(t)e^{y(t)}, \rho_{\ast}(t)]\subseteq [\rho_{min}e^{y(0)}, \rho_{max}]$ and $\overline{M}_{min}=\min\limits_{x\in  [\rho_{min}e^{y(0)}, \rho_{max}]}\left\{(-\frac{ \tanh \rho}{\rho})'\right\}>0$. 
Therefore,
\begin{align}
\frac{(1-e^{y(t)})'}{(1-e^{y(t)})}\le -\mu\Phi_{\ast}  \overline{M}_{min}\rho_{min}e^{y(0)}\label{4.15122}.
\end{align}
Integrating \re{4.15122} over $[0,t]$, we obtain
\begin{align}
1-e^{y(t)}\le (1-e^{y(0)})e^{-\mu\Phi_{\ast} \overline{M}_{min}\rho_{min} e^{y(0)}t}\qquad for ~    t>0\label{4.1511133}.
\end{align}
Taking $\delta=\min\{\mu\Phi_{\ast} {M}_{min}\rho_{min} , \mu\Phi_{\ast}\overline{M}_{min} \rho_{min} e^{y(0)}\}$ and $C=|1-e^{y(0)}|$,   from \re{4.15111} and \re{4.1511133}, we get \re{1.17}.

Finally, we prove the  uniqueness of  solution $\rho_{\ast}(t)$. Otherwise, using \re{1.17}, we have

$$
\left|\rho_{*}^{1}(t)-\rho_{*}^{2}(t)\right|\leqslant\left|\rho(t)-\rho_{*}^{1}(t)\right|+\left|\rho(t)-\rho_{*}^{2}(t)\right| \rightarrow 0 \qquad t\rightarrow \infty,
$$
i.e., $\rho_{*}^{1}(t)=\rho_{*}^{2}(t)$. Then we get the result.\hfill$\square$

\bibliographystyle{elsarticle-num}
\bibliography{3}

\end{document}